\documentclass[12pt,leqno]{amsart}

\usepackage{amscd}
\usepackage{amsmath}
\usepackage{amsfonts}
\usepackage{amssymb}
\usepackage{amsthm}
\usepackage{url}

\newtheorem*{defn}{Definition}

\title{Noether's problem for abelian extensions of cyclic $p$-groups II}
\author{Ivo M. Michailov}
\address{Faculty of Mathematics and Informatics, Shumen University "Episkop Konstantin Preslavski", Universitetska str. 115, 9700 Shumen, Bulgaria}
\email{ivo\_michailov@yahoo.com}
\date{\today}
\keywords{Noether's problem, Rationality problem, Meta-abelian group
actions, Metacyclic $p$-groups} \subjclass[2000]{primary 14E08
14M20; secondary 13A50,12F12}
\thanks{This work is partially supported by a project No RD-08-241/12.03.2013 of Shumen University}
\begin{document}
\baselineskip 20pt
\begin{abstract}
Let $K$ be a field and $G$ be a finite group. Let $G$ act on the
rational function field $K(x(g):g\in G)$ by $K$ automorphisms
defined by $g\cdot x(h)=x(gh)$ for any $g,h\in G$. Denote by $K(G)$
the fixed field $K(x(g):g\in G)^G$. Noether's problem then asks
whether $K(G)$ is rational (i.e., purely transcendental) over $K$.
Let $p$ be any prime and let $G$ be a $p$-group of exponent $p^e$.
Assume also that {\rm (i)} char $K = p>0$, or {\rm (ii)} char $K \ne
p$ and $K$ contains a primitive $p^e$-th root of unity. In this
paper we prove that if $G$ is any $p$-group of nilpotency class $2$,
which has the ABC (Abelian-By-Cyclic) property, then $K(G)$ is
rational over $K$. We also prove the rationality of $K(G)$ over $K$
for two $3$-generator $p$-groups $G$ of arbitrary  nilpotency class.
\end{abstract}

\maketitle
\newcommand{\Gal}{{\rm Gal}}
\newcommand{\Ker}{{\rm Ker}}
\newcommand{\GL}{{\rm GL}}
\newcommand{\Br}{{\rm Br}}
\newcommand{\lcm}{{\rm lcm}}
\newcommand{\ord}{{\rm ord}}
\renewcommand{\thefootnote}{\fnsymbol{footnote}}
\numberwithin{equation}{section}

\section{Introduction}
\label{1}

Let $K$ be any field. A field extension $L$ of $K$ is called
rational over $K$ (or $K$-rational, for short) if $L\simeq
K(x_1,\ldots,x_n)$ over $K$ for some integer $n$, with
$x_1,\ldots,x_n$ algebraically independent over $K$. Now let $G$ be
a finite group. Let $G$ act on the rational function field
$K(x(g):g\in G)$ by $K$ automorphisms defined by $g\cdot x(h)=x(gh)$
for any $g,h\in G$. Denote by $K(G)$ the fixed field $K(x(g):g\in
G)^G$. {\it Noether's problem} then asks whether $K(G)$ is rational
over $K$. This is related to the inverse Galois problem, to the
existence of generic $G$-Galois extensions over $K$, and to the
existence of versal $G$-torsors over $K$-rational field extensions
\cite[33.1, p.86]{Sw,Sa1,GMS}. Noether's problem for abelian groups
was studied extensively by Swan, Voskresenskii, Endo, Miyata and
Lenstra, etc. The reader is referred to Swan's paper for a survey of
this problem \cite{Sw}. Fischer's Theorem is a starting point of
investigating Noether's problem for finite abelian groups in
general.

\newtheorem{t1.1}{Theorem}[section]
\begin{t1.1}\label{t1.1}
{\rm (Fischer }\cite[Theorem 6.1]{Sw}{\rm )} Let $G$ be a finite
abelian group of exponent $e$. Assume that {\rm (i)} either char $K
= 0$ or char $K > 0$ with char $K\nmid e$, and {\rm (ii)} $K$
contains a primitive $e$-th root of unity. Then $K(G)$ is rational
over $K$.
\end{t1.1}

The next stage is the investigation of Noether's problem for finite
meta-abelian groups, and in particular metacyclic $p$-groups. Recall
that any metacyclic $p$-group $G$ is generated by two elements
$\sigma$ and $\tau$ with relations
$\sigma^{p^a}=1,\tau^{p^b}=\sigma^{p^c}$ and
$\tau^{-1}\sigma\tau=\sigma^{\varepsilon+\delta p^r}$ where
$\varepsilon=1$ if $p$ is odd, $\varepsilon=\pm 1$ if $p=2$,
$\delta=0,1$ and $a,b,c,r\geq 0$ are subject to some restrictions.
For the the description of these restrictions see e.g. \cite[p.
564]{Ka1}. The following Theorem of Kang generalizes Fischer's
Theorem for the metacyclic $p$-groups.

\newtheorem{t1.2}[t1.1]{Theorem}
\begin{t1.2}\label{t1.2}
{\rm (Kang}\cite[Theorem 1.5]{Ka1}{\rm )} Let $G$ be a metacyclic
$p$-group with exponent $p^e$, and let $K$ be any field such that
{\rm (i)} char $K = p$, or {\rm (ii)} char $K \ne p$ and $K$
contains a primitive $p^e$-th root of unity. Then $K(G)$ is rational
over $K$.
\end{t1.2}

Other results of Noether's problem for $p$-groups the reader can
find in \cite{CK,HuK,Ka2}.

\begin{defn} We say that a group $G$ has the ABC (Abelian-By-Cyclic)
property if $G$ has a normal abelian subgroup $H$ such that the
quotient group $G/H$ is cyclic
\end{defn}

Noether's problem for $p$-groups with the ABC property is still not
solved entirely. However, there are number of partial results that
have been obtained recently. We list some of them.

\newtheorem{t1.23}[t1.1]{Theorem}
\begin{t1.23}\label{t1.23}
{\rm (Haeuslein }\cite{Ha}{\rm )} Let $K$ be a field and $G$ be a
finite group. Assume that (i) $G$ contains an abelian normal
subgroup $H$ so that $G/H$ is cyclic of prime order $p$, (ii)
$\mathbb Z[\zeta_p]$ is a unique factorization domain, and (iii)
$\zeta_{p^e}\in K$ where $e$ is the exponent of $G$. If $G\to
\GL(V)$ is any finite-dimensional linear representation of $G$ over
$K$, then $K(V)^G$ is rational over $K$.
\end{t1.23}

\newtheorem{t1.24}[t1.1]{Theorem}
\begin{t1.24}\label{t1.24}
{\rm (Hajja }\cite{Haj}{\rm )} Let $K$ be a field and $G$ be a
finite group. Assume that (i) $G$ contains an abelian normal
subgroup $H$ so that $G/H$ is cyclic of order $n$, (ii) $\mathbb
Z[\zeta_n]$ is a unique factorization domain, and (iii) $K$ is
algebraically closed with $char K = 0$. If $G\to \GL(V)$ is any
finite-dimensional linear representation of $G$ over $K$, then
$K(V)^G$ is rational over $K$.
\end{t1.24}

\newtheorem{t1.22}[t1.1]{Theorem}
\begin{t1.22}\label{t1.22}
{\rm (}\cite[Theorem 1.4]{Ka3}{\rm )} Let $K$ be a field and $G$ be
a finite group. Assume that {\rm (i)} $G$ contains an abelian normal
subgroup $H$ so that $G/H$ is cyclic of order $n$, {\rm (ii)}
$\mathbb Z[\zeta_n]$ is a unique factorization domain, and {\rm
(iii)} $\zeta_{e}\in K$ where $e$ is the exponent of $G$. If
$G\rightarrow \GL(V)$ is any finite-dimensional linear
representation of $G$ over $K$, then $K(V)^G$ is rational over $K$.
\end{t1.22}

Note that those integers $n$ for which $\mathbb Z[\zeta_n]$ is a
unique factorization domain are determined by Masley and Montgomery.

\newtheorem{t1.25}[t1.1]{Theorem}
\begin{t1.25}\label{t1.25}
{\rm (Masley and Montgomery }\cite{MM}{\rm )}  $\mathbb Z[\zeta_n]$
is a unique factorization domain if and only if $1\leq n\leq 22$, or
$n = 24, 25, 26, 27, 28, 30, 32, 33, 34, 35, 36, 38, 40, 42, 45,
48,$ $50, 54, 60, 66, 70, 84, 90$.
\end{t1.25}

Therefore, Theorem \ref{t1.23} holds only for primes $p$ such that
$1\leq p\leq 19$. In a recent paper \cite{Mi} we managed to show
that the this condition can be waived, under some additional
assumptions regarding the structure of the abelian subgroup $H$.

\newtheorem{t1.4}[t1.1]{Theorem}
\begin{t1.4}\label{t1.4}
{\rm (Michailov}\cite[Theorem 1.8]{Mi}{\rm )} Let $p$ be an odd
prime, let $G$ be a group of order $p^n$ for $n\geq 2$ with an
abelian subgroup $H$ of order $p^{n-1}$, and let $G$ be of exponent
$p^e$. Choose any $\alpha\in G$ such that $\alpha$ generates $G/H$,
i.e., $\alpha\notin H,\alpha^p\in H$. Denote $H(p)=\{h\in H:
h^p=1,h\notin H^p\}$, and assume that $[H(p),\alpha]\subset H(p)$.
Denote by $G_{(i)}=[G,G_{(i-1)}]$ the lower central series for
$i\geq 1$ and $G_{(0)}=G$. Let the $p$-th lower central subgroup
$G_{(p)}$ be trivial. Assume that {\rm (i)} char $K = p>0$, or {\rm
(ii)} char $K \ne p$ and $K$ contains a primitive $p^e$-th root of
unity. Then $K(G)$ is rational over $K$.
\end{t1.4}

\newtheorem{t1.5}[t1.1]{Theorem}
\begin{t1.5}\label{t1.5}
{\rm (Michailov}\cite[Theorem 1.9]{Mi}{\rm )} Let $p$ be an odd
prime and let $G$ be a group of order $p^n$ for $n\leq 6$ which has
the ABC property. Let $G$ be of exponent $p^e$. Assume that {\rm
(i)} char $K = p>0$, or {\rm (ii)} char $K \ne p$ and $K$ contains a
primitive $p^e$-th root of unity. Then $K(G)$ is rational over $K$.
\end{t1.5}

The main purpose of this paper is to prove a generalization of the
latter results for all $p$-groups of nilpotency class $2$, which
have the ABC property. However, we should not ''over-generalize''
Theorem \ref{t1.5} to the case of any meta-abelian group because of
the following Theorem of Saltman.

\newtheorem{t1.3}[t1.1]{Theorem}
\begin{t1.3}\label{t1.3}
{\rm (Saltman }\cite{Sa2}{\rm )} For any prime number $p$ and for
any field $K$ with char $K \ne p$ (in particular, $K$ may be an
algebraically closed field), there is a meta-abelian $p$-group $G$
of order $p^9$ such that $K(G)$ is not rational over $K$.
\end{t1.3}

The main result of this paper is the following.

\newtheorem{main1}[t1.1]{Theorem}
\begin{main1}\label{main1}
For any prime $p$ let $G$ be a $p$-group of nilpotency class $2$,
which has the ABC property. Denote by $p^e$ the exponent of $G$.
Assume that {\rm (i)} char $K = p>0$, or {\rm (ii)} char $K \ne p$
and $K$ contains a primitive $p^e$-th root of unity. Then $K(G)$ is
rational over $K$.
\end{main1}

When $G$ has a nilpotency class bigger than $2$ it is much harder to
determine Noether's problem. In some cases, however, it is possible
to apply similar linearization techniques. We will do this for two
specific $p$-groups with three generators. Note that our method can
be applied for other $p$-groups with three generators.

Let $C_{p^a}$ be a cyclic group of order $p^a$ generated by the
element $\alpha$, and let $H$ be an abelian group generated by two
elements $\beta$ and $\gamma$ having orders $p^b$ and $p^c$,
respectively ($a,b,c\geq 1$). We are going to consider two groups
$G_1$ and $G_2$ which are in the middle of the exact sequence $1\to
H\to G_i\to C_{p^a}\to 1$. In particular, for the two groups we have
the relations
$\alpha^{p^a}\in\langle\beta,\gamma\rangle,\beta^{p^b}=\gamma^{p^c}=1,[\beta,\gamma]=1$.

Define
\begin{eqnarray*}
&&G_1=\langle\alpha,\beta,\gamma\mid
[\beta,\alpha]=1,~ [\gamma,\alpha]=\beta^x\gamma^{p^s},~ x\in\mathbb Z,~ s\in\mathbb N\rangle,\\
&&G_2=\langle\alpha,\beta,\gamma\mid [\beta,\alpha]=\beta^{p^r},~
[\gamma,\alpha]=\beta^x,~ x\in\mathbb Z,~ r\in\mathbb N\rangle,
\end{eqnarray*}

The second result of this paper is the following.

\newtheorem{main2}[t1.1]{Theorem}
\begin{main2}\label{main2}
For any prime $p$ let $G$ be isomorphic either to the group $G_1$ or
to $G_2$. Denote by $p^e$ the exponent of $G$. Assume that {\rm (i)}
char $K = p>0$, or {\rm (ii)} char $K \ne p$ and $K$ contains a
primitive $p^e$-th root of unity. Then $K(G)$ is rational over $K$.
\end{main2}

The key idea to prove Theorems \ref{main1} and \ref{main2} is to
find a faithful $G$-subspace $W$ of the regular representation space
$\bigoplus_{g\in G} K\cdot x(g)$ and to show that $W^G$ is rational
over $K$. The subspace $W$ is obtained as an induced representation
from $H$. By applying various linearizing techniques we then reduce
the rationality problem to another rationality problem which is
related either to cyclic or to metacyclic actions. The latter
actions can be linearized by using Lemma \ref{l2.7} or Theorem
\ref{t2.6}.

We organize this paper as follows. We recall some preliminaries in
Section \ref{2}. The proofs of Theorems \ref{main1} and \ref{main2}
are given respectively in Sections \ref{3} and \ref{4}.

\section{Preliminaries}
\label{2}

We list several results which will be used in the sequel.

\newtheorem{t2.1}{Theorem}[section]
\begin{t2.1}\label{t2.1}
{\rm (}\cite[Theorem 1]{HK}{\rm )} Let $G$ be a finite group acting
on $L(x_1,\dots,x_m)$, the rational function field of $m$ variables
over a field $L$ such that
\begin{description}
    \item [(i)] for any $\sigma\in G, \sigma(L)\subset L;$
    \item [(ii)] the restriction of the action of $G$ to $L$ is
    faithful;
    \item [(iii)] for any $\sigma\in G$,
    \begin{equation*}
\begin{pmatrix}
\sigma(x_1)\\
\vdots\\
\sigma(x_m)\\
\end{pmatrix}
=A(\sigma)\begin{pmatrix}
x_1\\
\vdots\\
x_m\\
\end{pmatrix}
+B(\sigma)
\end{equation*}
where $A(\sigma)\in\GL_m(L)$ and $B(\sigma)$ is $m\times 1$ matrix
over $L$. Then there exist $z_1,\dots,z_m\in L(x_1,\dots,x_m)$ so
that $L(x_1,\dots,x_m)^G=L^G(z_1,\dots,z_m)$ and $\sigma(z_i)=z_i$
for any $\sigma\in G$, any $1\leq i\leq m$.
\end{description}
\end{t2.1}

\newtheorem{t2.2}[t2.1]{Theorem}
\begin{t2.2}\label{t2.2}
{\rm (}\cite[Theorem 3.1]{AHK}{\rm )} Let $G$ be a finite group
acting on $L(x)$, the rational function field of one variable over a
field $L$. Assume that, for any $\sigma\in G,\sigma(L)\subset L$ and
$\sigma(x)=a_\sigma x+b_\sigma$ for any $a_\sigma,b_\sigma\in L$
with $a_\sigma\ne 0$. Then $L(x)^G=L^G(z)$ for some $z\in L[x]$.
\end{t2.2}

\newtheorem{t2.3}[t2.1]{Theorem}
\begin{t2.3}\label{t2.3}
{\rm (}\cite[Theorem 1.7]{CK}{\rm )} If $char K=p>0$ and $\widetilde
G$ is a finite $p$-group, then $K(G)$ is rational over $K$.
\end{t2.3}

The following Lemma can be extracted from some proofs in
\cite{Ka2,HuK}.

\newtheorem{l2.7}[t2.1]{Lemma}
\begin{l2.7}\label{l2.7}
Let $\langle\tau\rangle$ be a cyclic group of order $n>1$, acting on
$K(v_1,\dots,v_{n-1})$, the rational function field of $n-1$
variables over a field $K$ such that
\begin{eqnarray*}
\tau&:&v_1\mapsto v_2\mapsto\cdots\mapsto v_{n-1}\mapsto (v_1\cdots
v_{n-1})^{-1}\mapsto v_1.
\end{eqnarray*}
If $K$ contains a primitive $n$-th root of unity $\xi$, then
$K(v_1,\dots,v_{n-1})=K(s_1,\dots,s_{n-1})$ where $\tau:s_i\mapsto
\xi^is_i$ for $1\leq i\leq n-1$.
\end{l2.7}
\begin{proof}
Define $w_0=1+v_1+v_1v_2+\cdots+v_1v_2\cdots
v_{n-1},w_1=(1/w_0)-1/n,w_{i+1}=(v_1v_2\cdots v_i/w_0)-1/n$ for
$1\leq i\leq n-1$. Thus $K(v_1,\dots,v_{n-1})=K(w_1,\dots,w_n)$ with
$w_1+w_2+\cdots+w_n=0$ and
\begin{eqnarray*}
\tau&:&w_1\mapsto w_2\mapsto\cdots\mapsto w_{n-1}\mapsto w_n\mapsto
w_1.
\end{eqnarray*}
Define $s_i=\sum_{1\leq j\leq n}\xi^{-ij}w_j$ for $1\leq i\leq n-1$.
Then $K(w_1,\dots,w_n)=K(s_1,\dots,s_{n-1})$ and $\tau:s_i\mapsto
\xi^is_i$ for $1\leq i\leq n-1$.
\end{proof}

Now, let $G$ be any metacyclic $p$-group generated by two elements
$\sigma$ and $\tau$ with relations
$\sigma^{p^a}=1,\tau^{p^b}=\sigma^{p^c}$ and
$\tau^{-1}\sigma\tau=\sigma^{\varepsilon+\delta p^r}$ where
$\varepsilon=1$ if $p$ is odd, $\varepsilon=\pm 1$ if $p=2$,
$\delta=0,1$ and $a,b,c,r\geq 0$ are subject to some restrictions.
For the the description of these restrictions see e.g. \cite[p.
564]{Ka1}.

\newtheorem{t2.6}[t2.1]{Theorem}
\begin{t2.6}\label{t2.6}
{\rm (Kang }\cite[Theorem 4.1]{Ka1}{\rm )} Let $p$ be a prime
number, $m,n$ and $r$ are positive integers, $k=1+p^r$ if $(p,r)\ne
(2,1)$ (resp. $k=-1+2^r$ with $r\geq 2$). Let $G$ be a split
metacyclic $p$-group of order $p^{m+n}$ and exponent $p^e$ defined
by $G=\langle\sigma,\tau:
\sigma^{p^m}=\tau^{p^n}=1,\tau^{-1}\sigma\tau=\sigma^k\rangle$. Let
$K$ be any field such that $char K\ne p$ and $K$ contains a
primitive $p^e$-th root of unity, and let $\zeta$ be a primitive
$p^m$-th root of unity. Then $K(x_0,x_1,\dots,x_{p^n-1})^G$ is
rational over $K$, where $G$ acts on $x_0,\dots,x_{p^n-1}$ by
\begin{eqnarray*}
\sigma&:&x_i\mapsto \zeta^{k^i}x_i,\\
\tau&:&x_0\mapsto x_1\mapsto\cdots\mapsto x_{p^n-1}\mapsto x_0.
\end{eqnarray*}
\end{t2.6}

\section{Proof of Theorem \ref{main1}}
\label{3}

Let $G$ be generated by an abelian normal subgroup $H$ and an
element $\alpha$ such that $\alpha^{p^a}\in H$. Assume that
$H=\langle\alpha_1,\dots,\alpha_s\mid\alpha_i^{p^{a_i}}=1,1\leq
i\leq s\rangle$. We divide the proof into several steps. We are
going now to find a faithful representation of $G$.

\emph{Step 1.} Let $V$ be a $K$-vector space whose dual space $V^*$
is defined as $V^*=\bigoplus_{g\in G}K\cdot x(g)$ where $G$ acts on
$V^*$ by $h\cdot x(g)=x(hg)$ for any $h,g\in G$. Thus
$K(V)^G=K(x(g):g\in G)^G=K(G)$.

Define $X_1,X_2,\dots,X_{s}\in V^*$ by
\begin{equation*}
X_j=\sum_{\ell_1,\dots,\ell_{s}}x\left(\prod_{i\ne
j}\alpha_i^{\ell_i}\right),\quad \text{for}\ 1\leq j\leq s.
\end{equation*}

Note that $\alpha_i\cdot X_j=X_j$ for $j\ne i$. Let
$\zeta_{p^{a_i}}\in K$ be a primitive $p^{a_i}$-th root of unity for
$1\leq i\leq s$. Define $Y_1,Y_2,\dots,Y_{s}\in V^*$ by

\begin{equation*}
Y_i=\sum_{m=0}^{p^{a_i}-1}\zeta_{p^{a_i}}^{-m}\alpha_i^m\cdot X_i
\end{equation*}
for $1\leq i\leq s$.

It follows that {\allowdisplaybreaks\begin{align*} \alpha_i\ :\
&Y_i\mapsto\zeta_{p^{a_i}} Y_i,~ Y_j\mapsto Y_j,\ \text{for}\ j\ne
i.
\end{align*}}
Thus $\bigoplus_{1\leq j\leq s}K\cdot Y_j$ is a faithful
representation space of the subgroup
$H=\langle\alpha_1,\dots,\alpha_s\rangle$. The induced subspase $W$
depends on the relations in $G$. It is well known that the case
$\alpha^{p^a}\in H$ can easily be reduced to the case
$\alpha^{p^a}=1$ (see e.g. \cite[Proof of Theorem 1.8, Step 2]{Mi}).

Define $x_{ji}=\alpha^i\cdot Y_j$ for $1\leq j\leq s,0\leq i\leq
p^a-1$. Recall that $G$ is of nilpotency class $2$, so
$[H,\alpha]\leq Z(G)$, i.e., $[\alpha_j,\alpha]=\gamma_j\in Z(G)$
for $1\leq j\leq s$. It is not hard to see that
$$\alpha^{-i}\alpha_j\alpha^i=\alpha_j\gamma_j^i,\quad \text{for}\ 1\leq j\leq s,1\leq i\leq p^a-1.$$
We can now write the decomposition of $\gamma_j$ in $H$:
\begin{equation}\label{e1}
\gamma_j=\prod_{i=1}^s\alpha_i^{\alpha_{ij}p^{r_{ij}}}\quad
\text{for}\ 1\leq j\leq s,r_{ij}\geq 0,\alpha_{ij}\in\mathbb
Z,\gcd(\alpha_{ij},p)=1.
\end{equation}
It follows that {\allowdisplaybreaks\begin{align*}\alpha_{j}\ :\
&x_{ji}\mapsto\zeta_{p^{a_j}}\zeta_{p^{a_j}}^{i\alpha_{jj}p^{r_{jj}}}
x_{ji},~
x_{mi}\mapsto\zeta_{p^{a_m}}^{i\alpha_{mj}p^{r_{mj}}} x_{mi},\ \text{for}\ m\ne j,\\
\alpha\ :\ &x_{j0}\mapsto x_{j1}\mapsto\cdots\mapsto
x_{jp^a-1}\mapsto x_{j0},
\end{align*}}
where $0\leq i\leq p^a-1$ and $1\leq j\leq s$.

\emph{Step 2.} In this step we will find somewhat simpler actions of
$H$, which resemble the diagonal actions. First of all, for abuse of
notation we can assume that $\alpha_{ij}=1$ for all $i,j$. (We can
easily adjust the substitutions that follow in such a way that this
method works for arbitrary $\alpha_{ij}$'s, but the notations will
become cumbersome. Moreover, for the next steps we need not to know
the exact choices of primitive roots of unity.)

Next, without loss of generality we can assume that
$a_1-r_{11}=\max\{a_m-r_{mj}\mid 1\leq m,j\leq s\}$. Then we can
replace the generators of $H$ with the new generators
$\alpha_1,\alpha_1^{-p^{r_{12}-r_{11}}}\alpha_2,\dots,\alpha_1^{-p^{r_{1s}-r_{11}}}\alpha_s$.
The actions of the new generators are of the following type
{\allowdisplaybreaks\begin{align*}\alpha_1^{-p^{r_{1j}-r_{11}}}\alpha_{j}\
:\ &x_{1i}\mapsto\zeta_{p^{b_{1j}}} x_{1i},~
x_{mi}\mapsto\zeta_{p^{b_{mj}}}^i x_{mi},\ \text{for}\ m\ne j,
\end{align*}}
for some $b_{tj}\geq 0$, where $1\leq t\leq s,2\leq j\leq s$ and
$0\leq i\leq p^a-1$. Now, for $2\leq m\leq s$ define
$y_{mi}=x_{mi}x_{1i}^{a_1-r_{11}-a_m+p^{r_{m1}}}$. Then we have
{\allowdisplaybreaks\begin{align*}\alpha_1\ :\
&x_{1i}\mapsto\zeta_{p^{a_1}}\zeta_{p^{a_1}}^{ip^{r_{11}}} x_{1i},~
y_{mi}\mapsto\zeta_{p^{c_{m1}}} y_{mi},\ \text{for}\ m\ne
1,\\
\alpha_1^{-p^{r_{1j}-r_{11}}}\alpha_{j}\ :\
&x_{1i}\mapsto\zeta_{p^{c_{1j}}} x_{1i},~
y_{mi}\mapsto\zeta_{p^{d_{mj}}}\zeta_{p^{c_{mj}}}^i y_{mi},\
\text{for}\ m\ne 1,
\end{align*}}
for some $c_{tj}\geq 0~ (1\leq t,j\leq s)$ and $d_{tj}\geq 0 ~
(2\leq t,j\leq s)$, where $0\leq i\leq p^a-1$.

We can apply the same process with each generator
$\alpha_1^{-p^{r_{1j}-r_{11}}}\alpha_j$ for $2\leq j\leq s$. For
example, we may assume that $c_{22}=\max\{c_{mj}\mid 2\leq m,j\leq
s\}$. Then we can proceed in the same way for the generator
$\alpha_1^{-p^{r_{12}-r_{11}}}\alpha_{2}$. Repeating this process
with each generator we will finally obtain generators
$\beta_1,\dots,\beta_s$ of $H$ which act on the function field
$K(y_{ji}:1\leq j\leq s,1\leq i\leq p^a-1)=K(x_{ji}:1\leq j\leq
s,1\leq i\leq p^a-1)$ in this way
{\allowdisplaybreaks\begin{align*}\beta_j\ :\
&y_{ji}\mapsto\zeta_{p^{b_j}}\zeta_{p^{b_{jj}}}^i y_{ji},~
y_{mi}\mapsto\zeta_{p^{b_{mj}}} y_{mi},\ \text{for}\ m\ne j,
\end{align*}}
for some $b_j\geq 0,b_{tj}\geq 0$, where $1\leq t, j\leq s$ and
$0\leq i\leq p^a-1$. Clearly, the action of $\alpha$ is not changed:
\begin{align*}
\alpha\ :\ &y_{j0}\mapsto y_{j1}\mapsto\cdots\mapsto
y_{jp^a-1}\mapsto y_{j0},
\end{align*}
where $0\leq i\leq p^a-1$ and $1\leq j\leq s$.

Observe that from $1=[\alpha_j,\alpha^{p^a}]=\gamma_j^{p^a}$ it
follows that $a_m-r_{mj}\leq a$ for all $m,j$ (see \eqref{e1}).
Since each primitive root $\zeta_{p^{b_{mj}}}$ is obtained via
multiplication of roots of the type $\zeta_{p^{a_m-r_{mj}}}$, we
deduce that $b_{tj}\leq a$.

Clearly, $W=\bigoplus_{j,i}K\cdot y_{ij}\subset V^*$ is the induced
$G$-subspace obtained from $V$. Thus, by Theorem \ref{t2.1} it
suffices to show that $W^G$ is rational over $K$.

\emph{Step 3.} Assume that $b_{jj}\geq 1$, where $1\leq j\leq t$ for
some $t:1\leq t\leq s$. Denote $A=b_{11}$ and $\xi=\zeta_{p^{a-A}}$,
a primitive $p^{a-A}$-th root of unity. For $0\leq \ell\leq
p^{a-A}-1$ and $0\leq k\leq p^{A}-1$ define
$$u_{\ell k}=y_{1k}+\xi^\ell y_{1,k+p^{A}}+(\xi^\ell)^2
y_{1,k+2p^{A}}+\cdots+(\xi^\ell)^{p^{a-A}-1}
y_{1,k+(p^{a-A}-1)p^{A}}.$$

It is not hard to see that this is a well defined non-singular
transformation, so $K(u_{\ell k}:0\leq \ell\leq p^{a-A}-1,0\leq
k\leq p^{A}-1)=K(y_{1i}:0\leq i\leq p^a-1)$. Then the actions of the
generators of $G$ on $K(u_{\ell k})$ are
{\allowdisplaybreaks\begin{align*} \beta_1\ :\ &u_{\ell
k}\mapsto\zeta_{p^{b_1}}\zeta_{p^{A}}^k u_{\ell k},\\
\beta_m\ :\ &u_{\ell k}\mapsto\zeta_{p^{b_{1m}}} u_{\ell k},\quad \
\text{for}\ m\ne 1,\\
\alpha\ :\ &u_{\ell 0}\mapsto u_{\ell 1}\mapsto\cdots\mapsto u_{\ell
p^{A}-1}\mapsto \xi^{-\ell}u_{\ell 0},
\end{align*}}
where $0\leq \ell\leq p^{a-A}-1$ and $0\leq k\leq p^{A}-1$.

For $1\leq i\leq p^{A}-1$ define $v_{0i}=u_{0i}/u_{0i-1}$. For
$1\leq \ell\leq p^{a-A}-1$ and $0\leq k\leq p^{A}-1$ define $v_{\ell
k}=u_{\ell k}/u_{0 k}$. Thus $K(u_{\ell k}:0\leq \ell\leq
p^{a-A}-1,0\leq k\leq p^{A}-1)=K(u_{00},v_{\ell
k}:(k,\ell)\ne(0,0))$, and for every $g\in G$
\begin{equation}\label{e2}
g\cdot u_{00}\in K(v_{\ell k}:(k,\ell)\ne(0,0))\cdot u_{00},
\end{equation}
while the subfield $K(v_{\ell k}:(k,\ell)\ne(0,0))$ is invariant by
the action of $G$, i.e., {\allowdisplaybreaks\begin{align*} \beta_1\
:\ &v_{0i}\mapsto\zeta_{p^{A}}v_{0i},~ v_{\ell k}\mapsto v_{\ell
k},\quad \
\text{for}\ 1\leq i\leq p^{A}-1,0\leq k\leq p^{A}-1,\ell\ne 0\\
\beta_m\ :\ &v_{\ell k}\mapsto v_{\ell k},\quad \ \text{for}\ m\ne
1,0\leq
\ell\leq p^{a-A}-1,0\leq k\leq p^{A}-1,(k,\ell)\ne(0,0),\\
\alpha\ :\ &v_{01}\mapsto v_{02}\mapsto\cdots\mapsto
v_{0p^{A}-1}\mapsto (v_{01}\cdots v_{0p^{A}-1})^{-1},\\
&v_{\ell 0}\mapsto v_{\ell 1}\mapsto\cdots\mapsto v_{\ell
p^{A}-1}\mapsto \xi^{-\ell}v_{\ell 0}, \quad \ \text{for}\ 1\leq
\ell\leq p^{a-A}-1.
\end{align*}}
\emph{Step 4.} We can apply the same type of transformations as in
Step 3 for each $K(y_{ji}:0\leq i\leq p^a-1)$, where $2\leq j\leq
t$. For $t+1\leq j\leq s$ define $w_{ji}=y_{ji}/y_{ji-1}$, where
$1\leq i\leq p^a-1$. Denote by $\mathcal X$ the set of all new
variables (together with the variables of the type \eqref{e2}) and
by  $\mathcal X_j$ the subset of all variables that are not
invariant under $\beta_j$ and are not of the type \eqref{e2}.
Observe that $\beta_1$ leaves invariant all new variables except
$v_{0i}$'s and except the type \eqref{e2}, i.e. $\mathcal
X_1=\{v_{0i}:1\leq i\leq p^{A}-1\}$.

Define $w_1=v_{01}^{p^A}$ and $w_i=v_{0i}/v_{0i-1}$ for $2\leq i\leq
p^A-1$. Then $K(y_{ji}:0\leq i\leq p^a-1,1\leq j\leq
s)^{\langle\beta_1\rangle}=K(\mathcal
X)^{\langle\beta_1\rangle}=K(w_i,\mathcal X\setminus \mathcal
X_1,1\leq i\leq p^A-1)$. Note that for the variables of of the type
\eqref{e2} we apply Theorem \ref{t2.2}. The action of $\alpha$ on
$\mathcal X_1$ is {\allowdisplaybreaks
\begin{align*}
\alpha\ :\ &w_1\mapsto w_1w_2^{p^A},\\
&w_2\mapsto w_3\mapsto\cdots\mapsto w_{p^A-1}\mapsto
(w_1w_2^{p^A-1}w_3^{p^A-2}\cdots w_{p^A-1}^2)^{-1}\mapsto\\
&\mapsto w_1w_2^{p^A-2}w_3^{p^A-3}\cdots w_{p^A-2}^2w_{p^A-1}\mapsto
w_2.
\end{align*}}
Define $z_1=w_2,z_i=\alpha^i\cdot w_2$ for $2\leq i\leq p^A-1$. Now
the action of $\alpha$ is {\allowdisplaybreaks
\begin{align*}
\alpha\ :\ &z_1\mapsto z_2\mapsto\cdots\mapsto z_{p^A-1}\mapsto
(z_1z_2\cdots z_{p^A-1})^{-1}.
\end{align*}}
Since $w_1=(z_{p^A-1}z_1^{p^A-1}z_2^{p^A-2}\cdots
z_{p^A-2}^2)^{-1}$, we get that
$K(w_1,\dots,w_{p^A-1})=K(z_1,\dots,$ $z_{p^A-1})$. From Lemma
\ref{l2.7} it follows that the action of $\alpha$ on
$K(z_1,\dots,z_{p^A-1})$ can be linearized. Clearly, the action of
$\alpha$ on $K(v_{\ell i}:1\leq \ell\leq p^{a-A}-1,0\leq i\leq
p^A-1)$ is also linear.

Similarly, we can linearize the action of $\alpha$ on $K(\mathcal
X)^{\langle\beta_j\rangle}$ for each $j:2\leq j\leq t$. For
$j:t+1\leq j\leq s$ the action of $\alpha$ on $K(w_{ji})$ is
\begin{align*}
\alpha\ :\ &w_{j1}\mapsto w_{j2}\mapsto\cdots\mapsto
w_{jp^a-1}\mapsto (w_{j1}w_{j2}\cdots w_{jp^a-1})^{-1},
\end{align*}
which again can be linearized according to Lemma \ref{l2.7}.
Therefore, we obtain a linear action of $\alpha$ on $K(\mathcal
X)^H$. We are done.

\section{Proof of Theorem \ref{main2}}
\label{4}

Let $V$ be a $K$-vector space whose dual space $V^*$ is defined as
$V^*=\bigoplus_{g\in G}K\cdot x(g)$ where $G$ acts on $V^*$ by
$h\cdot x(g)=x(hg)$ for any $h,g\in G$. Thus $K(V)^G=K(x(g):g\in
G)^G=K(G)$. The key idea is to find a faithful $G$-subspace $W$ of
$V^*$ and to show that $W^G$ is rational over $K$. The subspace $W$
is obtained as an induced representation from $H$.

Define $X_1,X_2\in V^*$ by
\begin{equation*}
X_1=\sum_{i=0}^{p^b-1}x(\beta^i),~
X_2=\sum_{i=0}^{p^c-1}x(\gamma^i).
\end{equation*}
Note that $\beta\cdot X_1=X_1$ and $\gamma\cdot X_2=X_2$.

Let $\zeta_{p^b}\in K$ be a primitive $p^b$-th root of unity, and
let $\zeta_{p^c}\in K$ be a primitive $p^c$-th root of unity. Define
$Y_1,Y_2\in V^*$ by

\begin{equation*}
Y_1=\sum_{i=0}^{p^c-1}\zeta_{p^c}^{-1}\gamma^i\cdot X_1,~
Y_2=\sum_{i=0}^{p^b-1}\zeta_{p^b}^{-1}\beta^i\cdot X_2.
\end{equation*}

It follows that {\allowdisplaybreaks \begin{align*}
\beta\ :\ &Y_1\mapsto Y_1,~ Y_2\mapsto\zeta_{p^b} Y_2,\\
\gamma\ :\ &Y_1\mapsto\zeta_{p^c} Y_1,~ Y_2\mapsto Y_2.
\end{align*}}
Thus $K\cdot Y_1+K\cdot Y_2$ is a representation space of the
subgroup $H$. The induced subspase $W$ depends on the relations in
$G$.  It is well known that the case $\alpha^{p^a}=\beta^f\gamma^h$
can easily be reduced to the case $\alpha^{p^a}=1$ (see e.g.
\cite[Proof of Theorem 1.8, Step 2]{Mi}). Henceforth we will
consider separately $G_1$ and $G_2$.

\emph{Case I. $G\simeq G_1$.}

Define $x_i=\alpha^i\cdot Y_1,y_i=\alpha^i\cdot Y_2$ for $0\leq
i\leq p^a-1$. Calculations show that
$\gamma\alpha^i=\alpha^i\gamma^{k(i)}\beta^{xl(i)}$ for
$k(i)=1+\binom{i}{1}p^s+\binom{i}{2}p^{2s}+\cdots+\binom{i}{i}p^{is}=k^i$,
where $k=1+p^s$, and
$l(i)=\binom{i}{1}+\binom{i}{2}p^{s}+\cdots+\binom{i}{i}p^{(i-1)s}$.
We have now {\allowdisplaybreaks \begin{align*}
\beta\ :\ &x_i\mapsto x_i,~ y_i\mapsto\zeta_{p^b}y_i,\\
\gamma\ :\ &x_i\mapsto\zeta_{p^c}^{k^i} x_i,~ y_i\mapsto\zeta_{p^b}^{xl(i)}y_i,\\
\alpha\ :\ &x_0\mapsto x_1\mapsto\cdots\mapsto x_{p^a-1}\mapsto x_0,\\
&y_0\mapsto y_1\mapsto\cdots\mapsto y_{p^a-1}\mapsto y_0.
\end{align*}}
for $0\leq i\leq p^a-1$.

For $1\leq i\leq p^a-1$, define $u_i=x_i/x_{i-1}$ and
$v_i=y_i/y_{i-1}$. Thus $K(x_i,y_i:0\leq i\leq
p^a-1)=K(x_0,y_0,u_i,v_i:1\leq i\leq p^a-1)$ and for every $g\in G$
\begin{equation*}
g\cdot x_0\in K(u_i,v_i:1\leq i\leq p^a-1)\cdot x_0,~ g\cdot y_0\in
K(u_i,v_i:1\leq i\leq p^a-1)\cdot y_0,
\end{equation*}
while the subfield $K(u_i,v_i:1\leq i\leq p^a-1)$ is invariant by
the action of $G$. Thus $K(x_i,y_i:0\leq i\leq
p^a-1)^{G}=K(u_i,v_i:1\leq i\leq p^a-1)^{G}(u,v)$ for some $u,v$
such that $\alpha(v)=\beta(v)=\gamma(v)=v$ and
$\alpha(u)=\beta(u)=\gamma(u)=u$. Notice that $l(i)-l(i-1)=k^{i-1}$.
We have now {\allowdisplaybreaks
\begin{align}\label{e4.1}
\nonumber \beta\ :\ &u_i\mapsto u_i,~ v_i\mapsto v_i,\\
\gamma\ :\ &u_i\mapsto\zeta_{p^c}^{p^sk^{i-1}} u_i,~ v_i\mapsto\zeta_{p^b}^{xk^{i-1}} v_i,\\
\nonumber \alpha\ :\ &u_1\mapsto u_2\mapsto\cdots\mapsto u_{p^a-1}\mapsto (u_1u_2\cdots u_{p^a-1})^{-1},\\
\nonumber &v_1\mapsto v_2\mapsto\cdots\mapsto v_{p^a-1}\mapsto
(v_1v_2\cdots v_{p^a-1})^{-1},
\end{align}}
for $1\leq i\leq p^a-1$. From Theorem \ref{t2.2} it follows that if
$K(u_i,v_i:1\leq i\leq p^a-1)^{G}$ is rational over $K$, so is
$K(x_i,y_i:0\leq i\leq p^a-1)^{G}$ over $K$.

We can always write $x$ in the form $x=yp^t$ for $y\in\mathbb Z:
\gcd(y,p)=1$ and $t\geq 0$. Assume that $b-t\leq c-s$. For $1\leq
i\leq p^a-1$, define $w_i=v_i/u_i^{yp^{c-s-b+t}}$. Then
$K(u_i,v_i:1\leq i\leq p^a-1)=K(u_i,w_i:1\leq i\leq p^a-1)$ and
$\gamma(w_i)=w_i$ for $1\leq i\leq p^a-1$. (The other case $b-t>c-s$
is identical, since we can define $w_i=u_i/v_i^{zp^{b-t-c+s}}$ for
some $z$ such that $zy\equiv 1\pmod{p^{c-s}}$.)

Now, consider the metacyclic $p$-group $\widetilde
G=\langle\sigma,\tau:\sigma^{p^{c}}=\tau^{p^a}=1,\tau^{-1}\sigma\tau=\sigma^{k},k=1+p^s\rangle$.

Define $X=\sum_{0\leq j\leq
p^{c}-1}\zeta_{p^{c}}^{-j}x(\sigma^j),V_i=\tau^i X$ for $0\leq i\leq
p^a-1$. It follows that
\begin{eqnarray*}
\sigma&:&V_i\mapsto \zeta_{p^{c}}^{k^i}V_i,\\
\tau&:&V_0\mapsto V_1\mapsto\cdots\mapsto V_{p^a-1}\mapsto V_0.
\end{eqnarray*}
Note that $K(V_0,V_1,\dots,V_{p^a-1})^{\widetilde G}$ is rational by
Theorem \ref{t2.6}.

Define $U_i=V_i/V_{i-1}$ for $1\leq i\leq p^a-1$. Then
$K(V_0,V_1,\dots,V_{p^a-1})^{\widetilde G}=K(U_1,U_2,\dots,$
$U_{p^a-1})^{\widetilde G}(U)$ where
\begin{eqnarray*}
\sigma&:&U_i\mapsto \zeta_{p^{c}}^{k^i-k^{i-1}}U_i,~ U\mapsto U\\
\tau&:&U_1\mapsto U_2\mapsto\cdots\mapsto U_{p^a-1}\mapsto
(U_1U_2\cdots U_{p^a-1})^{-1},~ U\mapsto U.
\end{eqnarray*}

Notice that
$\zeta_{p^{c}}^{k^i-k^{i-1}}=\zeta_{p^{c}}^{p^sk^{i-1}}$. Compare
Formula \ref{e4.1} (i.e., the actions of $\gamma$ and $\alpha$ on
$K(u_i:1\leq i\leq p^a-1)$) with the actions of $\widetilde G$ on
$K(U_i:1\leq i\leq p^a-1)$. They are the same. Hence, according to
Theorem \ref{t2.6}, we get that $K(u_1,\dots,u_{p^a-1})^{G}(u)\cong
K(U_1,\dots,U_{p^a-1})^{\widetilde
G}(U)=K(V_0,V_1,\dots,V_{p^a-1})^{\widetilde G}$ is rational over
$K$. Since by Lemma \ref{l2.7} we can linearize the action of
$\alpha$ on $K(w_i:1\leq i\leq p^a-1)$, we obtain finally that
$K(u_i,w_i:1\leq i\leq p^a-1)^{\langle\gamma,\alpha\rangle}$ is
rational over $K$.

\emph{Case II. $G\simeq G_2$.}

Define $x_i=\alpha^i\cdot Y_1,y_i=\alpha^i\cdot Y_2$ for $0\leq
i\leq p^a-1$. Calculations show that
$\gamma\alpha^i=\alpha^i\gamma\beta^{xl(i)}$ for
$l(i)=\binom{i}{1}+\binom{i}{2}p^{r}+\cdots+\binom{i}{i}p^{(i-1)r}$,
and $\beta\alpha^i=\alpha^i\beta^{k(i)}$ for
$k(i)=1+\binom{i}{1}p^r+\binom{i}{2}p^{2r}+\cdots+\binom{i}{i}p^{ir}=k^i$,
where $k=1+p^r$. We have now {\allowdisplaybreaks
\begin{align*}
\beta\ :\ &x_i\mapsto x_i,~ y_i\mapsto\zeta_{p^b}^{k^i}y_i,\\
\gamma\ :\ &x_i\mapsto\zeta_{p^c} x_i,~ y_i\mapsto\zeta_{p^b}^{xl(i)}y_i,\\
\alpha\ :\ &x_0\mapsto x_1\mapsto\cdots\mapsto x_{p^a-1}\mapsto x_0,\\
&y_0\mapsto y_1\mapsto\cdots\mapsto y_{p^a-1}\mapsto y_0.
\end{align*}}
for $0\leq i\leq p^a-1$.

For $1\leq i\leq p^a-1$, define $u_i=x_i/x_{i-1}$ and
$v_i=y_i/y_{i-1}$. Thus $K(x_i,y_i:0\leq i\leq
p^a-1)=K(x_0,y_0,u_i,v_i:1\leq i\leq p^a-1)$ and for every $g\in G$
\begin{equation*}
g\cdot x_0\in K(u_i,v_i:1\leq i\leq p^a-1)\cdot x_0,~ g\cdot y_0\in
K(u_i,v_i:1\leq i\leq p^a-1)\cdot y_0,
\end{equation*}
while the subfield $K(u_i,v_i:1\leq i\leq p^a-1)$ is invariant by
the action of $G$. Thus $K(x_i,y_i:0\leq i\leq
p^a-1)^{G}=K(u_i,v_i:1\leq i\leq p^a-1)^{G}(u,v)$ for some $u,v$
such that $\alpha(v)=\beta(v)=\gamma(v)=v$ and
$\alpha(u)=\beta(u)=\gamma(u)=u$. Notice that $l(i)-l(i-1)=k^{i-1}$.
We have now {\allowdisplaybreaks
\begin{align*}
\beta\ :\ &u_i\mapsto u_i,~ v_i\mapsto\zeta_{p^b}^{p^rk^{i-1}} v_i,\\
\gamma\ :\ &u_i\mapsto u_i,~ v_i\mapsto\zeta_{p^b}^{xk^{i-1}} v_i,\\
\alpha\ :\ &u_1\mapsto u_2\mapsto\cdots\mapsto u_{p^a-1}\mapsto (u_1u_2\cdots u_{p^a-1})^{-1},\\
&v_1\mapsto v_2\mapsto\cdots\mapsto v_{p^a-1}\mapsto (v_1v_2\cdots
v_{p^a-1})^{-1},
\end{align*}}
for $1\leq i\leq p^a-1$. From Theorem \ref{t2.2} it follows that if
$K(u_i,v_i:1\leq i\leq p^a-1)^{G}$ is rational over $K$, so is
$K(x_i,y_i:0\leq i\leq p^a-1)^{G}$ over $K$.

We can always write $x$ in the form $x=yp^t$ for $y\in\mathbb Z:
\gcd(y,p)=1$ and $t\geq 0$. Assume that $r\leq t$. Since $\gamma$
acts in the same way as $\beta^{yp^{t-r}}$ on $K(u_i,v_i:1\leq i\leq
p^a-1)$, we find that $K(u_i,v_i:1\leq i\leq
p^a-1)^{G}=K(u_i,v_i:1\leq i\leq
p^2-1)^{\langle\beta,\alpha\rangle}$. (The other case $r>t$ is
identical, since $\beta$ acts in the same way as $\gamma^{zp^{r-t}}$
or some $z$ such that $zy\equiv 1\pmod{p^{b-r}}$.) The proof
henceforth is the same as Case I.

\end{document}